\documentclass[11pt]{paper}

\usepackage{amssymb,amsmath,enumerate,theorem,epsfig,rotating,graphics,changebar,eepic}
\usepackage{amsfonts}
\usepackage{a4,latexsym,parskip}
\usepackage{hyperref}
\hypersetup{pdfauthor=MS} \hypersetup{pdftitle=Fields}

\newcommand{\PP}{\mathbb{P}}
\newcommand{\A}{\mathbb{A}}

\newcommand{\C} [1][]{\mathbb{C}^{#1}}
\newcommand{\Q} [1] []{\mathbb{Q}_{#1}}
\newcommand{\N} [1][] {\mathbb{N}_{#1}}

\newcommand{\Z}{\mathbb{Z}}

\newcommand{\OO}{\mathcal{O}}

\newcommand{\qed}{\hfill \ensuremath{\Box}}

\newtheorem{Spezial-Theorem}{Theorem}[section]
\newtheorem{Spezial-Proposition}{Proposition}[section]
\theoremstyle{break} \newtheorem{Theorem}{Theorem}
\newtheorem{Proposition}[Theorem]{Proposition}
\newtheorem{Lemma}[Theorem]{Lemma}
\newtheorem{Example}[Theorem]{Example}

\newtheorem{Definition}[Theorem]{Definition}
\newtheorem{Corollary}[Theorem]{Corollary}

\newtheorem{Remark}[Theorem]{Remark}
\newtheorem{Question}[Theorem]{Question}

\begin{document}
\setlength{\unitlength}{1cm}

\title{Fields of definition of singular K3 surfaces}

\author{Matthias Sch\"utt}


\date{\today}
\maketitle

\abstract{This paper gives upper and lower bounds for the degree of the field of definition of a singular K3 surface, generalising a recent result by Shimada. We use work of Shioda-Mitani and Shioda-Inose and classical theory of complex multiplication.}




\section{Introduction}

Singular K3 surfaces have been paid considerable attention in algebraic and arithmetic geometry since the groundbreaking work of Shioda and Inose \cite{SI}. It follows from their theory that in many respects, singular K3 surfaces behave like elliptic curves with complex multiplication. For instance, they are defined over number fields, and some finiteness result applies (cf.~Thm.~\ref{Thm:Shafa}).

This paper is concerned with fields of definition of singular K3 surfaces. On the one hand, given an isomorphism class, we determine a model over a certain number field in Prop.~\ref{Prop:model}. This approach follows the ideas of Shioda and Inose in \cite{Inose}, \cite{Shioda} and \cite{SI}.

On the other hand, we derive a lower bound on the degree of the field of definition in Thm.~\ref{Thm:genus}. This result generalises the recent analysis of a special, yet important case by Shimada \cite[Thm.~3 (T)]{Shimada}. Our techniques employ class field theory (Sect.~\ref{s:CM}) and work of Shioda and Mitani  \cite{SM} on singular abelian surfaces (Sect.~\ref{s:abelian}). In some cases, Thm.~\ref{Thm:genus} can be combined with Lem.~\ref{Lem:conj} to strengthen the result (cf.~Sect.~\ref{s:Rem}).

We conclude the paper with some applications and remarks as to singular K3 surfaces over $\Q$ and  future directions in this subject.
Throughout the paper we will only be concerned with smooth complex projective varieties. Notation will be introduced wherever it is needed. A systematic summary can be found at the beginning of Sect.~\ref{s:abelian}.

\section{CM elliptic curves}
\label{s:CM}

In this section, we shall recall the classical theory of elliptic curves with CM. This will serve both as a motivation and as a toolbox for the similarly behaved singular K3 surfaces. The main reference for this section is the book of Shimura \cite{Shimura}.

Let $E$ denote an elliptic curve. It is well-known that $E$ can be defined over $\Q(j(E))$, and obviously this field is optimal. With respect to the group structure, $\Z\subseteq$ End$(E)$. An elliptic curve is said to have \emph{complex multiplication (CM)} if $\Z\subsetneq$ End$(E)$. In this case, End$(E)$ is an order $\OO$ in some imaginary quadratic field $K$. 

To $\OO$, there is associated the ring class field $H(\OO)=K(j(\OO))$ which is an abelian extension of $K$. Let $d_K$ denote the discriminant of $K$ and $f$ the conductor of $\OO$. Set $d=f^2d_K$. Consider the class group $Cl(d)$, consisting of positive-definite integral matrices of the form
\begin{eqnarray}\label{eq:Q}
Q=\begin{pmatrix}
2a & b\\
b & 2c
\end{pmatrix}
\end{eqnarray}
with discriminant $d=b^2-4ac$ up to the standard action of SL$_2(\Z)$. Dirichlet composition gives $Cl(d)$ the structure of an abelian group. By class theory, the Galois group of the extension $H(\OO)/K$ is isomorphic to $Cl(d)$. On the other hand, we have a natural map
\[
\Xi: \;\;Cl(d) \longrightarrow \{\text{elliptic curves with CM by }\OO\}/\cong
\]
sending $Q$ as in (\ref{eq:Q}) to the complex torus
\begin{eqnarray}
E_\tau=E_{\Z+\tau\Z}=\C/(\Z+\tau\Z)
\end{eqnarray}
with 
\begin{eqnarray}\label{eq:tau0}
\tau=\frac{-b+\sqrt{d}}{2a}.
\end{eqnarray}

\begin{Theorem}[cf.~{\cite[Thm.~7.7]{C}, \cite[Thm.~5.7]{Shimura}}]
The map $\Xi$ is a bijection. Any elliptic curve with CM by $\OO$ can be defined over $H(\OO)$. The representatives arising from $\Xi$ form a complete system of Galois conjugates over $K$. $\Xi$ is compatible with the Galois action, as given by the isomorphism Gal$(H(\OO)/K)\cong Cl(d)$.
\end{Theorem}

We shall also employ ad\'elic techniques. These will provide a uniform description of the Galois action on all elliptic curves with CM by an order in the fixed imaginary quadratic field $K$.

Let $\A_K^*$ denote the id\'ele group of $K$ and $K^\text{ab}$ the abelian closure of $K$. In the standard notation, there is a canonical surjection
\begin{eqnarray*}
[\;\cdot\;,K]: \A_K^* & \to & \text{Gal}(K^\text{ab}/K)\\
s\;\; & \mapsto & \;\;\;\;[s,K].
\end{eqnarray*}
Furthermore, the id\'eles carry a natural operation on the set of $\Z$-lattices in $K$ which is exhibited primewise by multiplication (cf.~\cite[\S 5.3]{Shimura}). For $s\in\A_K^*$ and a lattice $\Lambda$, we shall simply write $\Lambda\mapsto s\Lambda$. In terms of the class group $Cl(d)$ and the corresponding lattices $\Z+\tau\Z$, this action can also be interpreted as multiplication by $s\OO$. In this context, two $\Z$-lattices in $K$ are multiplied by multiplying generators, and they are identified if they agree up to homothety, i.e.~up to scaling in $K$.

\begin{Theorem}[Main theorem of CM {\cite[Thm.~5.4]{Shimura}}]\label{Thm:CM}
Let $E=E_\Lambda$ be a complex elliptic curve with CM by an order in $K$. Let $\sigma\in$ Aut$(\C/K)$ and $s\in\A_K^*$ such that $\sigma=[s,K]$ on $K^\text{ab}$. Then
\[
E^\sigma \cong E_{s^{-1}\Lambda}.
\]
\end{Theorem}

\section{Singular K3 surfaces}

In this section, we shall review the theory of singular K3 surfaces. Most of it goes back to the classical treatment by Shioda and Inose \cite{SI}.

\begin{Definition}
A surface $X$ is called \emph{singular}, if its Picard number attains the Lefschetz bound:
\[
\rho(X)=h^{1,1}(X).
\]\
\end{Definition}

In other words, a K3 surface $X$ is singular if $\rho(X)=20$.
One particular aspect of singular K3 surfaces is that they involve no
moduli. In fact, the general moduli space of complex K3 surfaces
with Picard number $\rho\geq\rho_0$ has dimension $20-\rho_0$.
Nevertheless, the singular K3 surfaces are everywhere dense in the
period domain of K3 surfaces (with respect to the analytic
topology).

\begin{Definition}
For a surface $X$, we define the \emph{transcendental lattice} $T_X$ as the orthogonal
complement of the N\'eron-Severi lattice $NS(X)$:
\[
T_X=NS(X)^{\bot}\subset H^2(X,\Z).
\]
\end{Definition}

For a K3 surface $X$, lattice theory predicts that $T_X$ is even with signature $(2,20-\rho(X))$. Hence, if $X$ is moreover singular, then $T_X$ is a positive definite even oriented lattice of rank two (with the orientation coming from the holomorphic $2$-form on $X$). One of Shioda-Inose's main results was the Torelli theorem for singular K3 surfaces:

\begin{Theorem}[Shioda-Inose {\cite[Thm.~4]{SI}}]\label{Thm:SI}
The map
\[
X \mapsto T_X
\]
induces a (1:1)-correspondence between isomophism classes of singular K3 surfaces and isomorphism classes of positive definite even oriented lattices of rank two.
\end{Theorem}

Shioda and Inose established the surjectivity of the above map by means of an explicit construction which is nowadays often referred to as \emph{Shioda-Inose structure}. Since this will be crucial to our arguments, we recall their ideas in the following.

A positive definite even oriented lattice $T$ of rank two is given by an intersection form $Q$ as in (\ref{eq:Q}) up to SL$_2(\Z)$. If $T=T_X$ for some surface $X$, we shall refer to its discriminant $d$ also as the discriminant of $X$. As before, we let $K=\Q(\sqrt{d})$ denote the imaginary quadratic field associated to $Q$ (or $X$). 

In a preceeding work \cite{SM}, Shioda and Mitani constructed a singular abelian surface $A$ for any given intersection form $Q$. This arose as product of two elliptic curves $E_{\tau_1}, E_{\tau_2}$ with
\begin{eqnarray}\label{eq:tau}
\tau_1=\tau=\frac{-b+\sqrt{d}}{2a} \;\;\text{ and }\;\; \tau_2=\frac{b+\sqrt{d}}{2}.
\end{eqnarray}
By construction, $E_{\tau_1}, E_{\tau_2}$ are isogenous elliptic curves with CM in $K$.

\begin{Theorem}[Shioda-Mitani {\cite[\S 3]{SM}}]\label{Thm:SM}
For $Q$ as in (\ref{eq:Q}), let $\tau_1, \tau_2$ as above. Then $A=E_{\tau_1}\times E_{\tau_2}$ is a singular abelian surface with intersection form $Q$ on $T_A$.
\end{Theorem}

For an abelian surface $A$, we can consider its Kummer surface Km$(A)$ which is K3. This results in multiplying the intersection form on the transcendental lattice by $2$:
\begin{eqnarray*}
T_{\text{Km} (A)} = T_A(2).
\end{eqnarray*}
In the singular case, i.e.~starting with $A=E_{\tau_1}\times E_{\tau_2}$, Shioda-Inose proceeded as follows to achieve the original intersection form: An elaborate investigation of the double Kummer pencil on Km$(A)$ produced a certain elliptic fibration. Then a quadratic base change gave rise to a singular K3 surface $X$ with the original intersection form (cf.~Sect.~\ref{s:model} for an explicit equation). We can rephrase this in terms of the corresponding deck transformation on $X$. This is a \emph{Nikulin involution}, i.e.~it has eight fixed points and leaves the holomorphic $2$-form invariant:

\begin{Theorem}[Shioda-Inose {\cite[Thm.~2]{SI}}]
Any singular K3 surface $X$ admits a Nikulin involution whose quotient has a Kummer surface Km$(A)$ as minimal resolution. In particular, $X$ is equipped with a Shioda-Inose structure
\[
\begin{array}{ccccc}
X &&&& A\\
& \searrow && \swarrow &\\
&& \text{Km}(A) &&
\end{array}
\] 
\end{Theorem}

As a consequence of the above construction, singular K3 surfaces in many respects behave like elliptic curves with CM. For instance, the following is immediate:

\begin{Corollary}[Shioda-Inose {\cite[Thm.~6]{SI}}]
Any complex singular K3 surface can be defined over some number field.
\end{Corollary}

By lattice and class theoretic means, \v Safarevi\v c widened this analogy to the extent of the following finiteness result:

\begin{Theorem}[\v Safarevi\v c {\cite[Thm.~1]{Shafa}}]\label{Thm:Shafa}
Fix $n\in\N$. Up to isomorphism, there is only a finite number of
complex singular K3 surfaces which possess a
model defined over a number field of degree at most $n$ over
$\Q$.
\end{Theorem}

It is the minimal degree of the fields of definition for singular K3 surfaces which we will be ultimately aiming at. As a first step towards the classification of all singular K3 surfaces for fixed degree $n$, this paper derives upper and lower bounds for this degree.

\section{Model over $\Q(j(\tau_1), j(\tau_2))$}
\label{s:model}

For a complex elliptic curve $E$, we know that $E$ can be defined over $\Q(j(E))$. In this section, we shall derive an analog statement for singular K3 surfaces. To simplify its proof, we shall normalize the $j$-map such that $j(i)=1 \;(i^2=-1)$. This adjustment does not affect our results.

\begin{Proposition}\label{Prop:model}
Let $X$ be a singular K3 surface with intersection form $Q$. Let $\tau_1, \tau_2$ as in (\ref{eq:tau}). Then $X$ has a model over $\Q(j(\tau_1), j(\tau_2))$ ($\subseteq K(j(\tau_2))$).
\end{Proposition}

The proof builds on a result of Inose \cite{Inose}. For a singular K3 surface $X$ as in the proposition, Inose derives a defining equation as a (non-smooth) quartic in $\PP^3$. This quartic has coefficients in the field $\Q(\alpha, \beta)$ where 
\[
\alpha^3 = j(\tau_1) j(\tau_2) \;\;\;\text{and}\;\;\; \beta^2 = (1-j(\tau_1)) (1-j(\tau_2)).
\]
The quartic polynomial can be interpreted as the elliptic fibration arising from the quadratic base change in Shioda-Inose's construction. In \cite{Sandwich}, Shioda refers to it as \emph{Inose's pencil}:
\[
X:\;\; y^2 = x^3 -3\, \alpha\, t^4\, x + t^5\, (t^2-2\,\beta\, t+1).
\]
To prove the proposition, we only have to exhibit a twist of the above fibration which is defined over the claimed smaller field. If $\alpha\,\beta\neq 0$, this is achieved in terms of the variable change
\[
t\mapsto \beta\, t,\;\;\; x\mapsto \frac\beta\alpha\, x,\;\;\;\ y\mapsto \sqrt{\frac{\beta^3}{\alpha^3}} y.
\]
Writing $A=\alpha^3, B=\beta^2$, the transformation results in the Weierstrass equation
\begin{eqnarray}\label{eq:A,B}
X:\;\;\; y^2 = x^3 -3\,A\,B\, t^4\, x + A\,B\,t^5\, (B\,t^2-2\,B\, t+1). 
\end{eqnarray}
If $\alpha\,\beta=0$, it suffices to substitute the zero entries in the above transformation by $1$. This proves Proposition \ref{Prop:model}. \qed

\begin{Remark}
Shioda-Inose's construction stays valid for any product abelian surface $A=E\times E'$ (cf.~\cite{Shioda}). The Kummer surface of $A$ can be obtained from $X$ via the base change $t\mapsto t^2$:
\begin{eqnarray}\label{eq:Kummer}
\text{Km}(A):\;\;\; y^2 = x^3 -3\,A\,B\, t^4\, x + A\,B\,t^4\, (B\,t^4-2\,B\, t^2+1). 
\end{eqnarray}
\end{Remark}

\begin{Corollary}\label{Cor:Kummer}
Let $X$ be a singular K3 surface with intersection form $Q$. If $Q$ is $2$-divisible, then $X$ can be defined over $\Q(j(\tau_1), j(\frac{\tau_2}2))$.
\end{Corollary}

\emph{Proof:} If $Q$ is $2$-divisible, then $X=\text{Km}(A)$ for the singular abelian surface $A$ with intersection form $\frac 12 Q$. In view of (\ref{eq:tau}), this leads to the claimed coefficients which we insert in (\ref{eq:Kummer}). \qed

The fibrations (\ref{eq:A,B}), (\ref{eq:Kummer}) allow us to realize any singular K3 surface over some predicted number field. However, this field $L$ clearly need not be optimal: For instance, if the two $j$-invariants are conjugate by some quadratic Galois automorphism $\sigma$, then $A$ and $B$ are fixed by $\sigma$. This will in some cases decrease the degree (see Sect.~\ref{ss:-23}). In fact, for any $r\leq 4$ there are examples where $L$ has degree $2^r$ over $\Q$, but $X$ can be defined over $\Q$. We will briefly comment on this phenomenon in Sect.~\ref{ss:spec}.

We end this section with an illustration of the interplay between Kummer construction and Shioda-Inose structure, as indicated by Cor.~\ref{Cor:Kummer}:

\begin{Example}[Fermat quartic]
The perhaps most classical singular K3 surface is the Fermat quartic in $\PP^3$
\begin{eqnarray*}\label{eq:Fermat}
S:\;\;\;\; x_0^4+x_1^4+x_2^4+x_3^4=0.
\end{eqnarray*}
In \cite{PSS}, relying on a result later proven by Cassels \cite{Cassels}, Pjatecki\u\i -\v Sapiro and  \v Safarevi\v c show that the intersection form on $T_S$ is 
\[
Q=\begin{pmatrix} 8 & 0\\ 0 & 8
\end{pmatrix}.
\]
Our previous considerations provide us with two further ways to describe $S$. On the one hand, $S=$Km$(E_i\times E_{2i})$. Both elliptic curves are defined over $\Q$, so the same holds for this model of $S$. On the other hand, we can apply the Shioda-Inose construction to $A=E_i\times E_{4i}$. In this case, the latter factor is only defined over a quadratic extension of $\Q$. Hence it is a priori not clear how to descend the model given by the special case analog of (\ref{eq:A,B}) to $\Q$.
\end{Example}

\section{Lower bound for the degree of the field of definition}
\label{s:lower}

In this section, we state our main result, giving a lower bound on the degree of the field of definition $L$ of a singular K3 surface $X$. Note that the extension $\Q(j(\tau_1), j(\tau_2))/\Q$ from Prop.~\ref{Prop:model} is in general not Galois, but the composition with the CM-field $K$ always is by class field theory. Therefore we will rather work with $LK$ and compute the degree 
\[
l=[LK:K]\mid [L:\Q].
\]
The question of complex conjugation (in $\C$) will be briefly addressed in Sect.~\ref{ss:conj}.

Let $X$ be a singular K3 surface, defined over some number field $L$. Assume that $L$ contains the CM-field $K$ and is Galois over $K$. In studying the Galois conjugates $X^\sigma$, we are particularly concerned with the transcendental lattices $T_{X^\sigma}$. 

\begin{Lemma}\label{Lem:genus}
For any $\sigma\in$ Gal$(L/K)$, $T_{X^\sigma}$ lies in the same genus as $T_X$.
\end{Lemma}

\emph{Proof:} The N\'eron-Severi lattice is a geometric invariant. In particular, the discriminant forms of $NS(X)$ and $NS(X^\sigma)$ agree. By \cite[Cor.~1.9.4]{N}, the transcendental lattices $T_X, T_{X^\sigma}$ lie in the same genus. \qed

Our main result is the following theorem which was first established by Shimada \cite[Thm.~3 (T)]{Shimada} in the special case of fundamental discriminant $d$ (i.e.~$d=d_K$):

\begin{Theorem}\label{Thm:genus}
In the  above notation, the set of isomorphism classes of the transcendental lattices $T_{X^\sigma}$ equals the genus of $T_X$:
\[
\{[T_{X^\sigma}]; \sigma\in\text{Gal}(L/K)\} = \text{genus of }\, T_X.
\]
\end{Theorem}

We note that our proof for general transcendental lattices differs significantly from Shimada's special case. Before going into details, we observe the theorem's implications on the degree $l$. Denote the degree of primitivity of $T_X$ by $m$ and let $d'$ such that $d=m^2d'$. Let $h$ be the class number of $T_X$, i.e.~$h=h(d')=\# Cl(d')$. Denote by $g$ the number of genera in $Cl(d')$. Then the number $n$ of classes per genus in $Cl(d')$ is given by
\[
n=\frac hg = \# Cl^2(d').
\]

\begin{Corollary}
The number of classes per genus divides $[L:K]$. In the above notation, $n\mid l$.
\end{Corollary}

It follows from Shioda-Inose's construction (cf.~Thm.~\ref{Thm:SI}) that in order to prove Thm.~\ref{Thm:genus}, it suffices to derive the analogous statement for the corresponding abelian surface $A=E_{\tau_1}\times E_{\tau_2}$:

\begin{Theorem}\label{Thm:abelian}
Let $A$ be a singular abelian surface over some number field. Let $K=\Q(\sqrt{d})$ be the CM-field of $A$. Then
\[
\{[T_{A^\sigma}]; \sigma\in\text{Aut}(\C/K)\} = \text{genus of } \,T_A.
\]
\end{Theorem}

Proving this theorem has the advantage that we can work with the explicit factors given by Thm.~\ref{Thm:SM}. Moreover, the Galois action on these elliptic curves is completely understood in terms of class field theory as reviewed in Sect.~\ref{s:CM}.

\section{Singular abelian surfaces}
\label{s:abelian}

In order to prove Thm.~\ref{Thm:abelian}, we shall recall further results from Shioda-Mitani's paper \cite{SM}. Beforehand, we recall and fix some notation related to the construction of Thm.~\ref{Thm:SM}. Let $Q$ be an even positive-definite binary quadratic form as in (\ref{eq:Q}) with discriminant $d$ and $\tau_1, \tau_2$ as in (\ref{eq:tau}):

\begin{tabular}{llll}
\textbf{Notation:} &  & $K=\Q(\sqrt{d})$ & imaginary quadratic field associated to $Q$\\
& & $d_K$ & discriminant of $K$\\
& & $m$ & degree of primitivity of $Q$\\
& & $Q'=\frac 1m Q$ & primitive even quadratic form associated to $Q$\\
& & $d'=\frac d{m^2}$ & discriminant of $Q'$\\
& & $Cl(Q)$ & class group of $Q$ (= forms of $Cl(d')$ multiplied by $m$)\\
& & $f$ & conductor of $Q$: $d=f^2d_K$\\
& & $f'$ & conductor of $Q'$: $d'=f'^2d_K$\\
& & $\OO$ &  order in $K$ of conductor $f$\\
& & $\OO'$ & order in $K$ of conductor $f'$.\\
& & $Cl(\OO')$ & class group of $\OO'$
\end{tabular}

It is immediate that upon varying $Q$ within its class group, $\tau_2$ does not change essentially,  i.e.~always $\Z+\tau_2\Z=\OO$. In other words, the isomorphism class of $A$ (or $T_A$) is completely reflected in $\tau_1$ (plus $d$), as it varies with $Q$, or equivalently, with $Q'$.

To state Shioda-Mitani's result, recall the notion of multiplication and equivalence of $\Z$-lattices in $K$ defined in Sect.~\ref{s:CM}. Note that the product of two lattices which correspond to quadratic forms in $Cl(d')$ and $Cl(d)$, corresponds to a quadratic form in $Cl(d')$.

\begin{Proposition}[Shioda-Mitani {\cite[Prop.~4.5]{SM}}]\label{Prop:abelian}
Let $A$ be the abelian surface given by $Q$ in Thm.~\ref{Thm:SM}. Let $\Lambda_1, \Lambda_2$ be $\Z$-lattices in $K$ of conductor $f_1, f_2$. Then
\[
A\cong E_{\Lambda_1}\times E_{\Lambda_2} \Leftrightarrow \begin{cases}
\text{(i)} & \Lambda_1\Lambda_2\sim \Z+\tau_1\Z,\\
\text{(ii)} & f_1f_2=mf'^2\,(=ff').\\
\end{cases}
\]
\end{Proposition}

(Shioda-Mitani's original formulation involved another condition on the right-hand side, but this is implied by the first equivalence.)

\emph{Proof of Thm.~\ref{Thm:abelian}:} Upon conjugating $A$, the resulting lattices of its factors stay in the same class group, so the second condition of Prop.~\ref{Prop:abelian} is trivially fulfilled. Alternatively, this follows from the invariance of the genus of $T_{A^\sigma}$ which is completely analogous to Lem.~\ref{Lem:genus}. Hence we only have to consider the first condition of Prop.~\ref{Prop:abelian}.

Let $\sigma\in$ Aut$(\C/K)$ and $s\in\A_K^*$ such that $\sigma=[s,K]$ on $K^\text{ab}$. Denote the lattice $\Lambda_j=\Z+\tau_j\Z$. By Thm.~\ref{Thm:CM},
\[
A^\sigma\cong E_{s^{-1}\Lambda_1}\times E_{s^{-1}\Lambda_2}.
\]
Then Prop.~\ref{Prop:abelian} and commutativity give the isomorphism
\begin{eqnarray}\label{eq:sigma}
A^\sigma\cong E_{s^{-2}\Lambda_1}\times E_{\Lambda_2}.
\end{eqnarray}
With $\sigma$ resp.~$s$ varying, the lattices $(s^{-1}\OO')^2$ cover the whole principal genus in the class group $Cl(\OO')$, since this consists exactly of the squares $Cl^2(\OO')$. Interpreting $s^{-2}\Lambda_1=(s^{-2}\OO')\Lambda_1=(s^{-1}\OO')^2\Lambda_1$, we derive that  $s^{-2}\Lambda_1$ always lies in the same genus. Since throughout this argument, we fix the lattice $\Lambda_2$ of the second factor in (\ref{eq:sigma}), this proves the claim of Thm.~\ref{Thm:abelian}. \qed

\section{Applications and Remarks}
\label{s:Rem}

\subsection{Complex conjugation}
\label{ss:conj}

In Sect.~\ref{s:lower}, we decided to consider the field of definition $L$ of a singular K3 surface only over the CM-field $K$ of the surface. This assured the extension to be Galois, but certainly at the loss of some information. Part of this information can be recovered by considering the complex conjugation $\iota$ of $\C$.

It is immediate that $(E_\tau)^\iota=E_{\bar\tau}$. This directly translates to a singular K3 surface $X$ (or singular abelian surface) with intersection form $Q$. We obtain that $X^\iota$ has intersection form $Q^{-1}$ with respect to the group law in its class group $Cl(Q)$. In terms of the coefficients as in (\ref{eq:Q}), this corresponds to $b\mapsto -b$.

\begin{Lemma}\label{Lem:conj}
Let $X$ be a singular K3 surface over some number field $L$ with intersection form $Q$. If $Q^2\neq 1$ in $Cl(Q)$, then $2|[L:\Q]$.
\end{Lemma}

The condition is most easily checked when $Q$ is in \emph{reduced form}, i.e.~in terms of (\ref{eq:Q}) $-a<b\leq a\leq c$ (and $b\geq 0$ if $a=c$). Then $Q$ is $2$-torsion in $Cl(Q)$, if and only if one of the three inequalities becomes equality by \cite[Lem.~3.10]{C}.

\subsection{Example: $d=-23$}
\label{ss:-23}

Let $d=-23$ and $K=\Q(\sqrt{d})$. Then $h(d)=3$ and actually
\[
Cl(d)=\left\{\begin{pmatrix}2 & 1\\1 & 12\end{pmatrix}, \begin{pmatrix}4 & -1\\-1 & 6\end{pmatrix}, \begin{pmatrix}4 & 1\\1 & 6\end{pmatrix}\right\}
\]
consists of only one genus. This is the first case where $\Q(j(\OO_K))/\Q$ is not Galois, since the conjugates of $j(\OO_K)$ are clearly not real. Inserting $\tau_1=\tau_2=\frac{-1+\sqrt{d}}2$ into (\ref{eq:A,B}), we find that the singular K3 surface with the first intersection form is defined over $\Q(j(\OO_K))$. However, combining Thm.~\ref{Thm:genus} and Lem.~\ref{Lem:conj}, we see that the other two singular K3 surfaces can only be defined over the Hilbert class field $H=K(j(\OO_K))$. 

This observation agrees with the following fact: On the one hand, there is some quadratic $\sigma\in$ Gal$(H/K)$ such that we can write these other surfaces as in (\ref{eq:A,B}) with $j(\tau_2)=j(\tau_1)^\sigma$. Nevertheless, the fixed field of $\sigma$ still is just $\Q$, so this does not decrease the degree of the field of definition.

Since $h(-92)=3$, we obtain exactly the same results for singular K3 surfaces with discriminant $-92$. For the imprimitive case, apply Cor.~\ref{Cor:Kummer}.

\subsection{Concluding remarks}
\label{ss:spec}

For a singular K3 surface, we can always combine Thm.~\ref{Thm:genus} and Lem.~\ref{Lem:conj} to obtain a lower bound on its field of definition. Sec.~\ref{ss:-23} gave three examples where this lower bound coincided with the upper bound from Prop.~\ref{Prop:model}. This can be easily generalised:

\begin{Lemma}\label{Lem:bounds}
Let $X$ be a singular K3 surface with intersection form $Q$, discriminant $d$ and degree of primitivity $m$. Let $\tau_1$ as in (\ref{eq:tau}). Assume that $(d,m)$ or $(\frac d4, \frac m2)$ is to be found in the following table where $r$ always denotes an odd integer:

\begin{center}
\begin{tabular}{|c|c|}
\hline
$d$ & $m$\\
\hline\hline
$-4$ & $1$\\
$-8$ &\\
$-16$ &\\
$-p^r, p\equiv 3\mod 4$ &\\
$-4p^r, p\equiv 3\mod 4$ &\\
\hline
$-12$ & $2$\\
$-16$ &\\
$-4p^r, p\equiv 7\mod 8$ &\\
\hline
$-27$ & $3$\\
\hline
\end{tabular}
\end{center}

If $Q=1$ in $Cl(Q)$, then $\Q(j(\tau_1))$ is the minimal field of definition of $X$. Otherwise the minimal field of definition is $K(j(\tau_1))$.
\end{Lemma}

\emph{Proof:} 
The lemma follows whenever the mentioned upper and lower bounds coincide. In the primitive case, this is equivalent to the condition that $Cl(d)$ consists only of one genus. By \cite[Thm.s 3.11, 3.15]{C}, this gives the first entries. Generally, the imprimitive case additionally requires $h(d)=h(d')$. Thus we obtain the other entries from \cite[Cor.~7.28]{C}. This proves the claim for the pairs $(d,m)$.

If $2|m$, we can furthermore employ the Kummer construction from Cor.~\ref{Cor:Kummer}. This directly gives the claim for pairs $(\frac d4, \frac m2)$. \qed

If there is more than one genus, one might still hope that the lower bound is always attained in the primitive case (and maybe even if $\frac 12 Q$ is primitive). So far, this hope is particularly supported by examples with discriminant of relatively small absolute value. Mostly, these examples are extremal elliptic surfaces (cf.~\cite{BM}, \cite{S-Rocky}). 

There has been particular interest in singular K3 surfaces over $\Q$. By Thm.~\ref{Thm:genus}, the class group in this case is only $2$-torsion. One easily finds 101 discriminants satisfying this condition. By \cite{Wb}, there is at most one further discriminant (of enormously large absolute value). 

The number of corresponding imaginary quadratic fields is 65. The problem of finding singular K3 surfaces over $\Q$ just for these fields (regardless of the actual discriminant) is related to the following modularity question which was formulated independently by Mazur and van Straten: 

\begin{Question}[Mazur, van Straten]
Is any newform of weight 3 with rational Fourier coefficients associated to a singular K3 surface over $\Q$?
\end{Question}
With this respect, the $2$-torsion statement has already been obtained in \cite[Prop.~13.1]{S-CM} where one can also find a list of the 65 fields and the newforms (up to twisting).

So far, we know singular K3 surfaces over $\Q$ for 43 of these fields. Remarkably, the class groups in these cases are $2$-torsion of rank up to $4$. In other words, the class number becomes as large as $16$. In contrast, there is only one apparent general way to descend the field of definition. This essentially consists in the Weil restriction from an elliptic curve to an abelian surface because then we can apply the Shioda-Inose construction to this surface (cf.~the discussion succeeding Cor.~\ref{Cor:Kummer} and in Sect.~\ref{ss:-23}). However, this technique a priori only decreases the rank of $2$-torsion in the Galois group by $1$.


A general solution to the problem of the field of definition of a singular K3 surface seems still far away (although Lem.~\ref{Lem:bounds} gives a complete answer for a good portion of singular K3 surfaces). The main difficulty seems to lie in the imprimitive case: Replacing $Q$ by a multiple $mQ$ gives
\[
\tau_1(mQ)=\tau_1(Q),\;\;\; \tau_2(mQ)=m \tau_2(Q).
\]
In particular, the class number stays unchanged. In other words, our lower bound from Thm.~\ref{Thm:genus} is fixed. Meanwhile the upper bound of Prop.~\ref{Prop:model} increases with $m$, since $\tau_2$ changes. On the other hand, Thm.~\ref{Thm:Shafa} tells us that the degree of the field of definition does in fact increase with the degree of primitivity. It is our intention to pursue this issue in a future work. 

We conclude the paper with the following remark which underlines the subtleties involved in the imprimitive case: There is a singular K3 surface over $\Q$ with degree of primitivity as large as $30$. 
Beukers and Montanus in \cite{BM} found an elliptic K3 surface with configuration $[3,3,4,4,5,5]$. This has intersection form
\[
Q=30\begin{pmatrix}2 & 0\\0 & 2\end{pmatrix}.
\]

\vspace{0.8cm}

\textbf{Acknowledgement:} The paper was written while the author enjoyed the kind hospitality of Harvard University. Generous funding from DFG under research grant Schu 2266/2-1 is gratefully acknowledged. The author would like to thank N.~Elkies, B.~van Geemen, B.~Hassett, K.~Hulek and T.~Shioda for many stimulating discussions and the anonymous referee for helpful comments and suggestions.

\vspace{0.8cm}

\vspace{0.8cm}

Matthias Sch\"utt\\
Mathematics Department\\
Harvard University\\
1 Oxford Street\\
Cambridge, MA 02138\\
USA\\
{\tt mschuett@math.harvard.edu}

\end{document}